\documentclass{article}
\usepackage[utf8]{inputenc}

\usepackage[a4paper, total={6in, 8in}]{geometry}
\usepackage{amsmath,amssymb,amsthm}
\usepackage{hyperref}
\usepackage{graphicx}
\usepackage{xcolor}
\usepackage{algorithm}
\usepackage{algpseudocode}
\usepackage{multirow}
\usepackage{authblk}
\usepackage{caption}
\usepackage{subcaption}

\usepackage{amsmath}
\usepackage{amssymb}
\usepackage{amsthm}
\usepackage{amssymb}
\usepackage{bbm}
\usepackage{hyperref}
\usepackage{color}

\newtheorem{theorem}{Theorem}
\newtheorem{lemma}{Lemma}
\newtheorem{remark}{Remark}
\def \ER {Erd\H{o}s-R\'enyi }

\newcommand{\ceil}[1]{\lceil #1 \rceil}
\newcommand{\mb}[1]{\mathbb{#1}}

\newcommand{\brac}[1]{\left(#1\right)}
\newcommand{\cbrac}[1]{\left\{#1\right\}}
\newcommand{\sbrac}[1]{\left[#1\right]}
\newcommand{\indic}[1]{\mathbbm{1}{\brac{#1}}}
\newcommand{\expect}[2][]{\mathbb{E}_{#1}\sbrac{#2}}
\newcommand{\prob}[2][]{\mathbb{P}_{#1}\brac{#2}}
\newcommand{\define}{\triangleq}
\newcommand{\abs}[1]{\lvert #1 \rvert}

\newcommand{\mf}[1]{\mathbf{#1}}
\providecommand{\keywords}[1]
{
  \small	
  \textbf{\textit{Keywords---}} #1
}

\title{Phase Transitions in Biased Opinion Dynamics with 2-choices Rule}
\author{Arpan Mukhopadhyay\\
Department of Computer Science\\
University of Warwick}
\date{}

\begin{document}

\maketitle

\begin{abstract}
    We consider a model of binary opinion dynamics where one opinion is inherently `superior' than the other and social agents exhibit a `bias' towards the superior alternative. Specifically, it is assumed that an agent updates its choice to the superior alternative with  probability $\alpha >0$ irrespective of its current opinion and the opinions of the other agents. With probability $1-\alpha$ it adopts the majority opinion among two randomly sampled neighbours and itself. We are interested in the time it takes for the network to converge to a consensus on the superior alternative. In a complete graph of size $n$, we show that irrespective of the initial configuration of the network, the average time to reach consensus scales as $\Theta(n \log n)$ when  the bias parameter $\alpha$ is sufficiently high, i.e., $\alpha > \alpha_c$ where $\alpha_c$ is a threshold parameter that is uniquely characterised. When the bias is low, i.e., when $\alpha \in (0,\alpha_c]$, we show that the same rate of convergence can only be achieved if the initial proportion of agents with the superior opinion is above certain threshold $p_c(\alpha)$. If this is not the case, then we show that the network takes $\Omega(\exp(\Theta(n)))$ time on average to reach consensus. 
\end{abstract}


\section{Introduction}
\label{sec:intro}

Opinion dynamical models are used extensively in statistical physics and computer science to study the effects of different local interaction rules on the adoption  of new technologies and products. 
One key question in this context is how fast can a new/superior technology replace an old/outdated technology in a network of connected of agents?  
Classical opinion dynamical models with two competing opinions (or technologies) assume the opinions to be indistinguishable; indeed under both voter and majority rule models, agents update their opinions purely based on the opinions of other agents in their neighbourhoods without exhibiting any preference for any opinion. 
However, to capture the inherent superiority of one opinion over another, 
% we need models in which agents prefer one opinion over another; 
we need to incorporate some form of ``bias'' towards one of the two opinions. 

Opinion dynamical models with bias have been studied recently in~\cite{arpan_JSP_20,mukhopadhyay2016binary,ijcai2020p8,lesfari2022biased}.
A strong form of bias is considered in~\cite{ijcai2020p8,lesfari2022biased}; a bias parameter $\alpha \in (0,1)$ is introduced; with probability $\alpha$ an agent is assumed to perform a {\em biased update} in which it adopts the superior opinion  independently of its current opinion and the opinions of all other agents; with probability $1-\alpha$ the agent performs a {\em regular update} in which it adopts the majority opinion among {\em all} agents in its neighbourhood. 
The later update rule is often referred to as the {\em majority rule} in the literature.
Clearly, in this model, any network will eventually reach consensus on the superior opinion. It is shown in~\cite{ijcai2020p8} that the expected time taken for the network to reach consensus on the superior opinion is exponential in the minimum degree of the underlying  graph. Thus, in graphs where the neighbourhood size of each agent is proportional to the network size, the consensus time grows exponentially with the network size. This naturally prompts the question if simpler rules exist under which consensus can be achieved faster.

To address this question, we consider a simpler rule
called the {\em 2-choices rule}~\cite{cooper_two_choices} for {\em regular updates}. Under the 2-choices rule, an agent samples two other agents from its neighbourhood and updates to the majority opinion among the sampled agents and the agent itself.
For graphs where the neighbourhood sizes are large, this modified rule greatly reduces the communication overhead associated with computing the majority opinion since an agent no longer needs to know the opinions of all other agents in its neighbourhood; only knowing the opinions of the two sampled agents suffices. 
% in the modified model it needs to be done for a group constant size an agent now only has to find the majority opinion among a constant number of agents rather than having to find the majority opinion in a large group whose size grows with the network size $n$.
The 2-choices rule also reduces the chance that an updating agent adopts the worse alternative when a majority of its neighbours have chosen this alternative. Hence, this modification should not only reduce the communication overhead, but also `facilitate' consensus on the superior alternative. In this paper, we analytically characterise the improvement to the speed of consensus brought about by this modified rule.

Specifically, we show that if the network is a complete graph and the bias parameter $\alpha$ is sufficiently high ($\alpha > 1/9$), then consensus is achieved on the superior opinion in $\Theta(n \log n)$ time, where $n$ denotes the number of agents in the network. When bias is small (i.e., when $\alpha < 1/9$) we show that the consensus time depends on the initial configuration of the network. More specifically, when the bias parameter $\alpha$ is small, fast consensus (i.e., $\Theta(n \log n)$ time) is  achieved only if the initial proportion $p$ of agents with the superior opinion is above a certain threshold value $p_c(\alpha)$ (explicitly characterised). If $p < p_c(\alpha)$ we show that the network takes exponentially long time to reach consensus. Thus, the speed with which the network reaches consensus on the superior opinion undergoes a sharp {\em phase transition} depending on the values of $\alpha$ and $p$. 
Through simulations, we observe similar behaviour on other classes of graphs, e.g., on random $d$-regular graphs both with $d=\Theta(\log n)$ and \ER graphs with edge probability $\Theta(\log n/n)$. To establish our theoretical results, we use a novel characterisation of the expected number of visits of a random walk to a given state using a branching process. We expect this technique to be useful in the analysis of other interacting particle system models.
Thus, in summary, our contributions are the following:

\begin{itemize}
    \item (Fast consensus) For complete graphs, we show the existence of a sharp threshold such that if the bias parameter $\alpha$ is above the threshold, then consensus is achieved on the superior opinion starting from any initial configuration in $\Theta(n\log n)$ time on average. If the bias parameter $\alpha$ is below the threshold, we show the existence of another sharp threshold such that if the initial proportion $p$ of agents with the superior opinion is above this threshold, then consensus can be achieved in $\Theta(n \log n)$ time.
    
    \item (Slow consensus) We show that when both the bias parameter $\alpha$ and the initial proportion $p$ of agents with the superior opinion are below their corresponding thresholds, the average consensus time on complete graphs is $\Omega(\exp(\Theta(n)))$, i.e., grows exponentially with the network size.
    
    \item (Other classes of graphs) Through extensive simulations we study consensus on other classes of graphs. Specifically, we observe similar behaviour on random $d$-regular graphs with $d=\Theta(\log n)$ and on \ER graphs with edge probability $\log n/n$.
    For $d$-regular graphs with constant degrees $d=O(1)$ we do observe a phase transition but the behaviour below criticality is different from that in complete graphs or other dense graphs studied in this paper.
\end{itemize}

\subsection{Related Literature}

The simplest model of opinion dynamics studied in the literature is the {\em voter model} where an agent simply copies the opinion of an agent sampled randomly from its neighbourhood. Thus, in the voter model, the probability that an agent adopts a specific opinion is equal to the proportion of agents having the same opinion in its neighbourhood. Due to the linearity of the resulting dynamics and its duality with coalescing random walks, the voter model has been extensively studied in the literature. The duality between the voter model and coalescing random walks was first observed independently in~\cite{holley1975ergodic} and~\cite{clifford1973model}. Using this duality, the model has been analysed on different classes of graphs such as regular lattices~\cite{cox1989coalescing}, random $d$-regular graphs~\cite{cooper2013coalescing}, and ER graphs~\cite{nakata1999probabilistic}. It is known that for connected graphs the probability of reaching consensus on a specific opinion is proportional to the initial volume (sum of degrees) of nodes having that opinion. It is also known that the mean consensus time for the asynchronous version of the voter model (where in each round only one randomly sampled agent updates its opinion) is $\Omega(n^2)$ for most graphs.

Another important model of opinion dynamics is the majority rule model wherein an agent adopts the majority opinion among all agents in its neighbourhood. It was shown in~\cite{mossel2014majority} that with high probability for a family of expander graphs with sufficiently large spectral gap, the majority dynamics leads to a consensus on the opinion having the initial majority, provided that the imbalance between the initial majority opinion and the alternate opinion is sufficiently high. Bounds on the consensus time for the majority rule model was obtained in~\cite{zehmakan2020opinion} for expanders and Erdos-Renyi random graphs. For expanders with sufficiently large spectral gaps, it was shown that consensus can be achieved on the initial majority opinion in $O(\log n)$ steps in the synchronous model where all agents update in each round. For ER graphs, it was shown that if the edge probability is above the connectivity threshold of $\log n/n$, then consensus can be achieved on the initial majority opinion in constant number of rounds. The majority rule model has also been studied for other classes of graphs such as finite lattices~\cite{schonmann1990finite}, random regular graphs~\cite{gartner2018majority}, and infinite trees~\cite{kanoria2011majority}.

Although the majority rule leads to faster consensus on many classes of graphs, it requires an agent to know the states of all other agents in its neighbourhood. This may be too computationally expensive when neighbourhood sizes are large. A simpler alternative is to consider the 2-choices rule where an agent only samples two random neighbours and changes to the majority opinion among the sampled agents and the agent itself. Rules similar to the 2-choices rule, where groups of agents are formed at each instant and all agents in the group update to the majority opinion within the group, were analysed in the physics literature~\cite{redner_two_choices,galam2002minority,chen2005consensus}. A generalisation of 2-choices rule, where the updating agent samples $m$ agents from its neighbourhood and only changes its opinion if $d$ or more of the sampled agent differ from the updating agent, was analysed in the continuous time in~\cite{cruise2014probabilistic}. The 2-choices rule that we consider in this paper was first analysed for random $d$-regular graphs and expanders in~\cite{cooper_two_choices}. It was shown that consensus can be achieved in $O(\log n)$ time with high probability on the initial majority opinion provided that the initial imbalance is sufficiently high.

Opinion dynamical models with bias have been considered in the recent literature. These models are designed to capture the superiority of one alternative over another. Accordingly, in these models, the agents exhibit some form of bias towards the superior opinion.
In~\cite{arpan_JSP_20,mukhopadhyay2016binary} a weak form of bias is considered. Here, agents  with the superior opinion update with a lesser frequency than agents with the alternative opinion.
It has been shown in~\cite{arpan_JSP_20} that, under the voter rule and the 2-choices rule,  consensus is achieved in $O(\log n)$ time for complete graphs.
For the voter model, the probability of achieving consensus on the superior opinion approaches to one as the network size grows. For the 2-choices rule, consensus is achieved on the superior opinion with high probability only when the initial proportion of agents with the superior opinion is above a certain threshold. 

The form bias studied in this paper 
is introduced recently in~\cite{ijcai2020p8,biased_info_sciences_2022}. These papers show that on dense graphs the speed of consensus can be slow (depending on the minimum degree of a node) if the agents follow the majority rule during a regular update. The papers also analyse the voter rule under this form bias and {show} that  consensus with the voter rule  can be achieved in $O(n \log n)$ time. 
{Our model is different from these models as we consider the 2-choices rule as the main update rule.  
Furthermore, we study the dynamics under the 2-choices rule as a function of the bias parameter $\alpha$ as well as the initial proportion $p$ of agents with the superior opinion.
 This is unlike the previous papers~\cite{ijcai2020p8, biased_info_sciences_2022}, where the dynamics is studied
as a function of $\alpha$ only with a fixed value of $p$ (specifically,  $p=0$).
%Furthermore, we study the dynamics under the 2-choices rule as a function of not only the bias parameter $\alpha$ and but also as a function of the initial proportion $p$ of agents having the superior opinion. 
We show that phase transitions can occur with respect to both $\alpha$ and $p$}. Similar phase transitions have been recently reported in~\cite{cruciani2021phase,d2022phase,cruciani2021phase1} for noisy versions of the $k$-majority dynamics. However, the focus of these papers is a state of `near consensus' where both opinions co-exist but the proportion of agents with one of the opinions is arbitrarily small. In contrast, our results focus on the state of full consensus in which the inferior opinion is completely eliminated. Furthermore, we obtain bounds on the mean of the consensus time that are stronger than high probability bounds obtained in previous papers. The technique used here is also quite different from those used in the earlier works; while earlier works use concentration around the mean drift, we use suitably constructed branching processes to obtain bounds on the number of visits to different states.

\subsection{Organisation}

The rest of the paper is organised as follows. In Section~\ref{sec:model}, we introduce the model studied in this paper. Next, in Section~\ref{sec:analysis} we present the theoretical analysis of the model for complete graphs. Section~\ref{sec:numerics} provides numerical results to support our theoretical findings and also demonstrates similar behaviour for other classes of graphs. Finally, we conclude the paper in Section~\ref{sec:conclusion}.

% However, in a network which is both large and dense, neighbourhood size of each agent is large and as a result it is often infeasible for an agent to determine the majority opinion in its neighbourhood by looking at the opinion every agent in its neighbourhood.

% This motivates us to consider a simpler update rule which is computationally less expensive.
% In particular, we consider 
% Motivated by the above, we consider a variation of the model considered in~\cite{biased_info_sciences_2022} where, instead of looking at all agents in its neighbourhood, an agent with probability $1-\alpha$ chooses two random agents in its neighbourhood and changes its opinion if both the sampled agents have the same opinion opposite to that of the agent.

\section{Model}
\label{sec:model}

In this section, we describe the model studied in this paper. 
The model consists of a network of $n$ agents described by
an undirected graph $G_n=(V_n,E_n)$, where
the nodes in $V_n$ (with $\abs{V_n}=n$) represent the  agents
and the edges in $E_n$ represent the connections between the agents.
For each agent (node) $u \in V_n$,
we denote by $N_u = \cbrac{v:(u,v)\in E_n}$,
the set of neighbours of $u$.
% We call the elements of $\partial u$
% the neighbours of $u$.

Time is assumed to be discrete and at each discrete time instant $t \in \mathbb{Z}_+$ each agent is assumed to
have an opinion in the set $\{0,1\}$.
Without loss of generality, we assume that $1$ is the superior opinion.  
Let $X_u(t) \in \cbrac{0,1}$ denote the opinion
of agent $u$ at time $t$.
At $t=0$, the opinions of the agents are initialised
such that $\abs{\cbrac{u\in V_n: X_u(0)=1}}=\ceil{pn}$ for some $p \in [0,1)$. Hence,
$p$ fraction of agents initially have the superior opinion $1$. 

At each instant $t \geq 0$, an agent, sampled uniformly at random,
updates its opinion:
with probability $\alpha$, it performs a {\em biased update}
in which it adopts the superior opinion irrespective of its current opinion and the opinions of all other agents;
with probability $1-\alpha$ it performs a {\em regular update} following the {\em 2-choices rule} in which the updating agent samples two neighbours uniformly at random (with replacement\footnote{
Note that sampling with or without replacement
does not make any difference when the neighbourhood
size is proportional to $n$ since
the probability of choosing the same neighbour twice
tends to zero as $n \to \infty$.}) and adopts the majority opinion among the sampled agents and the agent itself.
Therefore, if $U(t)$ denotes the randomly sampled agent at time $t$, then at time $t+1$ the opinion of the agent is given by

\begin{equation}
    \label{eq:update_rule}
    X_{U(t)}(t+1)=\begin{cases}
                        1 & \text{ w.p. } \alpha, \\
                        M(t) & \text{ w.p. }  1-\alpha,
                   \end{cases}
\end{equation}
where $M(t)=\indic{X_{U(t)}(t)+X_{N_1(t)}(t)+X_{N_2(t)}(t)\geq 2}$ denotes the majority opinion among two randomly sampled neighbours $N_1(t)$ and $N_2(t)$ of $U(t)$ and the agent $U(t)$ itself.
The parameter $\alpha \in (0,1]$ represents the bias towards the superior opinion and is referred to as the {\em bias parameter} of the model. 

% To perform an update, the sampled agent $S$ 
% further samples two of its neighbours uniformly at
% random (with replacement)\footnote{
% Note that sampling with or without replacement
% does not make any difference when the neighbourhood
% size is proportional to $N$ since
% the probability of choosing the same neighbour twice
% tends to zero as $N \to \infty$.} from the set $\partial S(t)$.
% Let $N_1(t)$ and $N_2(t)$ denote these sampled neighbours.
% Let $M(t)=\indic{X_{S(t)}(t)+X_{N_1(t)}(t)+X_{N_2(t)}(t)\geq 2}$ be the majority opinion among $S(t), N_1(t)$, and $N_2(t)$.
% The agent $S(t)$ updates its opinion according to 
% the following rule:

% \begin{equation}
%     X_{S(t)}(t+1)=\begin{cases}
%                         1 & \text{ w.p. } \alpha, \\
%                         M(t) & \text{ w.p. }  1-\alpha,
%                   \end{cases}
% \end{equation}
% %
% where $\alpha \in (0,1]$ is called the {\em bias probability}
% of the model. Clearly, an agent chooses opinion $1$ irrespective of its neighbourhood with probability $\alpha$
% and chooses the majority opinion among its sampled neighbours with probability $1-\alpha$.
% Under the above rule, it is easy to see that for any
% network $G_N$, the agents converge to a consensus on opinion $1$ starting from any state in a finite time with probability one.
% Hence, in this case, $1$ is called the dominant opinion of 
% the agents.

The state of the network at any time $t \geq 0$ can be represented by the vector $\mf{X}(t)=(X_u(t), u \in V_n)$ of opinions of all the agents. The process $\mf{X}(\cdot)$ is Markov on the state space $\cbrac{0,1}^n$ with a single absorbing state $\mf{1}$ where all agents have opinion $1$. We refer to this absorbing state as the {\em consensus state}. Since it is possible to reach the consensus state from any other state in a finite number of steps, with probability one the chain is absorbed in the consensus state in a finite time. We refer to this time as the {\em consensus time}. The objective of the rest of the paper is to analyse the mean consensus time for different values of the parameters
$n,\alpha,p$ and for different classes of graphs.
% We are interested in characterising
% the consensus time $T_N$ for different network topologies
% as a function of the network size $N$. In particular,
% we are interested in the scaling law of $T_N$ either
% in expectation or with high probability.

\section{Analysis for complete graphs}
\label{sec:analysis}

In this section, we assume that $G_n$ is a complete graph
on $n$ nodes. 
% We also assume for simplicity that during a regular update
% an agent can sample itself in addition to sampling its neighbours.
For complete graphs,
the process $\bar X^n=(\bar X^n(t)=\abs{\cbrac{u\in V_n: X_u(t)=1}}, t\geq 0)$ counting the number of agents
with opinion $1$ is a Markov chain on $\mb Z_+$
with absorbing state $n$.
For the chain $\bar X^n$, 
the transition probability $\tilde{p}_{i,j}$ from
state $i$ to $j$ is given by $\tilde p_{i,j}={p}_{i,j}+o(1/n)$
where

\begin{equation}
\label{eq:transition_prob}
    p_{i,j}=\begin{cases}
                \brac{1-\frac{i}{n}}\brac{\alpha+(1-\alpha)\brac{\frac{i}{n}}^2}, &\text{ if } j=i+1,\\
                \frac{i}{n}(1-\alpha)\brac{1-\frac{i}{n}}^2, &\text{ if } j=i-1,\\
                1-p_{i,i+1}-p_{i,i-1}, &\text{ if } j=i,\\
                0, &\text{ otherwise.}
            \end{cases}
\end{equation}
% \arpan{Explain how these transition probabilities are obtained}
Note that $p_{i,j}$ denotes the transition probability of a slightly modified chain $\bar Y^n$ (with the same absorbing state $n$) where an agent during a regular update can sample itself 
in addition to its neighbours. 
Below we analyse this modified chain as it is simpler to do so and {the asymptotic (in $n$) results we obtain for this modified chain also hold for the original chain (see Remark~\ref{rem:mod_chain} for more explanation).}

Below we characterise the scaling law
of $\bar T_n(p)=\expect[\ceil{pn}]{T_n}$, where $\expect[x]{\cdot}$
denotes expectation conditioned on $\bar Y^n(0)=x$
and $T_k=\inf\cbrac{t\geq 0: \bar Y^n(t)=k}$ denotes the first time the network reaches state $k$.
Our main result is described in Theorem~\ref{thm:complete}
below.

\begin{theorem} 
\label{thm:complete}
For complete graphs we have the following
\begin{enumerate}
    \item (Fast consensus) For each $\alpha \in (1/9,1)$, we have $\bar T_n(p) = \Theta(n \log n)$ for all $p \in [0,1)$.
    
    \item (Fast consensus) For each $\alpha \in (0, 1/9)$, there exists $p_c(\alpha) \in (0,1)$ such that if $p \in [p_c(\alpha),1)$, then $\bar T_n(p)=\Theta(n \log n)$.
    Furthermore, the threshold $p_c(\alpha)$ is the unique solution in the range $(\bar x_\alpha,1)$ of the equation

    \begin{equation}
        \label{eq:p_cutoff}
        \int_{\underbar{x}_\alpha}^{p_c(\alpha)} \log(f_\alpha(x))dx =0,
    \end{equation}
    where $f_\alpha(x)=\frac{(1-\alpha)x(1-x)}{\alpha+(1-\alpha)x^2}$,
    $\underbar{x}_{\alpha}=\frac{1}{4}\brac{1-\sqrt{1-\frac{8\alpha}{1-\alpha}}}$ and
    $\bar x_\alpha =\frac{1}{4}\brac{1+\sqrt{1-\frac{8\alpha}{1-\alpha}}}$.
    
    \item (Slow consensus) For each $\alpha \in (0,1/9)$ and $p \in [0,p_c(\alpha))$, we have $\bar T_n(p)=\Omega(\exp(\Theta(n)))$.
\end{enumerate}
\end{theorem}

The above theorem implies that when the bias parameter is sufficiently high ($\alpha >1/9$), the network quickly (in $\Theta(n \log n)$ time) reaches consensus on the superior opinion irrespective of its initial state. This is in sharp contrast to the result of~\cite{biased_info_sciences_2022} where the consensus time is exponential in $n$ for all values of the bias parameter $\alpha$. The theorem further implies that even with low value of the bias parameter $\alpha$ (for $\alpha < 1/9$) fast consensus can be achieved as long as the initial proportion $p$ of agents with the superior opinion is above a threshold value denoted by $p_c(\alpha)$. We explicitly characterise this threshold $p_c(\alpha)$ required to ensure fast consensus when the bias is low. If the bias parameter is low ($\alpha < 1/9$) and the initial proportion of agents with the superior opinion is below the threshold, i.e., $p < p_c(\alpha)$, then the mean consensus time is exponential in $n$ which corresponds to slow speed of convergence. Thus, our model exhibits rich behaviour in terms of the parameters $\alpha$ and $p$.

We breakdown the proof of the above theorem into several simpler steps. The first step is to make the following simple observation which expresses the mean consensus time as a function of the number of visits of the chain $\bar Y^n$ to different states in its state-space $\cbrac{0,1,\ldots,n}$.

\begin{lemma}
\label{lem:master}
Let $Z_k$ denote the number of visits of the chain $\bar Y^n$ to the state $k \in \cbrac{0,1,\ldots,n}$ before absorption. Then, for any starting state $x =\bar{X}^n(0)\in \cbrac{0,1,\ldots,n}$,
the average consensus time is given by

\begin{equation}
    \expect[x]{T_n} =\sum_{k=0}^{n-1} \frac{\expect[x]{Z_k}}{(1-p_{k,k})}, \label{eq:master_relation}
\end{equation}
where $p_{k,k}$ is the transition probability from state $k$ to itself given by~\eqref{eq:transition_prob}.
\end{lemma}

\begin{proof}
Observe that the consensus time $T_n$ can be written as

\begin{equation}
    \label{eq:T_sum}
    T_n=\sum_{k=0}^{n-1}\sum_{j=1}^{Z_k} M_{k,j},
\end{equation}
where $Z_k$ denotes the number of visits to state $k$
before absorption and $M_{k,j}$ denotes the 
time spent in state $k$ in the $j^{\textrm{th}}$ visit.
Clearly, the random variables $Z_k$ and $(M_{k,j})_{j\geq 1}$ are independent of
each other. Furthermore, $(M_{k,j})_{j\geq 1}$
is a sequence of i.i.d. random variables with 
geometric distribution given by
$\prob[x]{M_{k,j}=i}=p_{k,k}^{i-1}(1-p_{k,k})$ for all $j$.
Hence, applying Wald's identity to~\eqref{eq:T_sum} we have

\begin{align}
\bar T_{n}(p) &=\expect[x]{T_n} \nonumber\\
                    &=\sum_{k=0}^{n-1}\expect[x]{\sum_{j=1}^{Z_k}M_{k,j}}\nonumber\\
					&=\sum_{k=0}^{n-1} \expect[x]{Z_k} \expect[x]{M_{k,j}} \nonumber\\
				    &=\sum_{k=0}^{n-1} \frac{\expect[x]{Z_k}}{(1-p_{k,k})}, \label{eq:consensus}
\end{align}
where the last step follows from
the fact that $\expect[x]{M_{k,j}}=1/(1-p_{k,k})$
for each $j \in [Z_k]$.
\end{proof}

From the above lemma, it is evident that in order to obtain bounds on $\expect[x]{T_n}$ we need to obtain bounds the expected number of visits $\expect[x]{Z_k}$ to different states and the transition probabilities $p_{k,k}$. Using this approach, we first obtain a lower bound on $\bar T_n(p)$.

\begin{lemma}
\label{lem:lower_bound}
For any $\alpha \in (0,1)$ and $p \in (0,1)$ we have $\bar T_n(p)=\Omega(n \log n)$.
\end{lemma}

\begin{proof}
From~\eqref{eq:transition_prob} we obtain

\begin{align*}
    1-p_{k,k} &= \brac{1-\frac{k}{n}}\brac{\alpha+(1-\alpha)\frac{k}{n}}\\
              &\leq \max\brac{\alpha,1-\alpha} \brac{1-\frac{k}{n}}\brac{1+\frac{k}{n}}
\end{align*}
Furthermore, we have $\expect[x]{Z_k} \geq \indic{k \geq x}$
since states $k \geq x$ are visited at least
once. Hence, using the above two inequalities
in~\eqref{eq:master_relation} we obtain

\begin{align*}
    \bar T_n(p) & \geq \sum_{k=0}^{n-1} \frac{\indic{k \geq \ceil{np}}}{\max\brac{\alpha,1-\alpha} \brac{1-\frac{k}{n}}\brac{1+\frac{k}{n}}}\\
    & = \frac{n}{2\max\brac{\alpha,1-\alpha}}\sum_{k=\ceil{np}}^{n-1}\brac{\frac{1}{n-k}+\frac{1}{n+k}}\\
    & > \frac{n}{2\max\brac{\alpha,1-\alpha}}\sum_{k=\ceil{np}}^{n-1}\brac{\frac{1}{n-k}}\\
    & \geq  \frac{n}{2\max\brac{\alpha,1-\alpha}}\log\brac{n-\ceil{np}+1},
\end{align*}
which completes the proof.
\end{proof}

To obtain an upper bound on the mean consensus time $\bar T_n(p)$ similarly, we need an upper bound on the expected number of visits $\expect[x]{Z_k}$ to state $k$ for each $k \in \cbrac{0,1,2\ldots,n}$. We obtain such an upper bound using a technique developed in~\cite{arpan_JSP_20} where the number of visits to each state is expressed as a function of a branching process. Specifically, we define $\zeta_k$ to be the number of jumps from state $k$
to state $k-1$ for any $k\in \cbrac{1,2,\ldots,n}$. Clearly, $\zeta_n=0$ and, if the starting state is $x=\bar Y^n(0) \in [1,n)$, then $\zeta_k$ satisfies the following recursion:

\begin{equation}
\label{eq:zeta}
    \zeta_k=\begin{cases}
              \sum_{l=0}^{\zeta_{k+1}}\xi_{l,k}, &\text{ if } k \in \{x,x+1,\ldots,n-1\},\\
              \sum_{l=1}^{\zeta_{k+1}}\xi_{l,k}, &\text{ if } k \in \{0,1,\ldots,x-1\},
            \end{cases}
\end{equation}
where $\xi_{l,k}$ denotes the number of jumps from state $k$
to state $k-1$ between the $l^{\textrm{th}}$ and $(l+1)^{\textrm{th}}$ visit to state $k+1$.
Note that the first sum in~\eqref{eq:zeta} (for $k \geq x$)
starts from $l=0$ whereas the second some starts from
$l=1$.
This is because for $k\geq x$, jumps from $k$ to $k-1$ can occur even before
the first jump from $k+1$ to $k$. 
This is not the case for states $k < x$. Furthermore,
$\cbrac{\xi_{l,k}}_{l\geq 0}$ is a sequence of i.i.d.
random variables independent of $\zeta_{k+1}$.
Hence, the sequence $\{\xi_{n-k}\}_k$ defines
a branching process. Moreover, the number of visits to state $k$
can be written as 

\begin{equation}
    \label{eq:z_k}
    {Z_k} = \begin{cases}
                         1+\zeta_k+\zeta_{k+1}, &\text{ if } k \in \{x,x+1,\ldots,n-1\},\\
                         \zeta_k+\zeta_{k+1}, &\text{ if } k \in \{0,1,\ldots,x-1\},
                      \end{cases}
\end{equation}
where the additional one appears in the first case (for $k \geq x$)
because the states $k \geq x$ are visited at least once before the chain is absorbed in state $n$. It is the above characterisation of the $Z_k$ that we shall use to obtain
upper bounds on $\expect[x]{Z_k}$. Precisely, we obtain bounds on $\expect[x]{\zeta_k}$
to bound $\expect[x]{Z_k}$.
In the lemma below, we express $\expect[x]{\zeta_k}$ in terms of the transition probabilities of the chain $\bar{Y}^n$.

\begin{lemma}
\label{lem:zeta_k}
For any $k\in \cbrac{1,2,\ldots,n}$, let $\zeta_k$ denote the number of jumps of the chain $\bar Y^n$ from state $k$ to state $k-1$. We have $\zeta_n=0$ and 

\begin{equation}
\label{eq:solved_expect_zeta}
\expect[x]{\zeta_{k}}=\begin{cases}
			\sum_{t=k}^{n-1} \prod_{i=k}^{t} \frac{p_{i,i-1}}{p_{i,i+1}}, \quad \text{for } x \leq k \leq n-1,\\
			\brac{\prod_{i=k}^{x-1} \frac{p_{i,i-1}}{p_{i,i+1}}}\expect[x]{\zeta_{x}}, \quad \text{for } 0 \leq k < x,
			\end{cases}
\end{equation}
where for each $i, j \in \{0,1,\ldots,n\}$, $p_{i,j}$ denotes the transition probability of the Markov chain $\bar Y^n$ from state $i$ to state $j$ and is given by~\eqref{eq:transition_prob}.
\end{lemma}

\begin{proof}
We observe that
$\zeta_{k+1}$ and $\cbrac{\xi_{l,k}}_{l\geq 0}$ are independent of each other and
$\cbrac{\xi_{l,k}}_{l\geq 0}$ is a sequence of i.i.d.
random variables with $$\prob[x]{\xi_{l,k}=i}=\brac{\frac{p_{k,k-1}}{1-p_{k,k}}}^i\brac{\frac{p_{k,k+1}}{1-p_{k,k}}}$$ for all $i \in \mb{Z}_+$ and all $l \geq 0$.
Hence, we have 
$$\expect[x]{\xi_{l,k}}=\frac{p_{k,k-1}}{p_{k,k+1}},$$ for all $l \geq 0$.
Applying Wald's identity to~\eqref{eq:zeta}
we obtain the following recursions

\begin{equation}
\label{eq:expect_zeta}
    \expect[x]{\zeta_k}=\begin{cases}
              \brac{1+\expect[x]{\zeta_{k+1}}}\frac{p_{k,k-1}}{p_{k,k+1}}, &\text{ if } k \geq x,\\
              \expect[x]{\zeta_{k+1}}\frac{p_{k,k-1}}{p_{k,k+1}}, &\text{ if } k < x,
            \end{cases}
\end{equation}
Upon solving the above recursions with boundary condition
${\zeta_n}=0$, we obtain desired result.
\end{proof}

From~\eqref{eq:solved_expect_zeta}, we observe that the ratio $p_{i,i-1}/p_{i,i+1}$ plays a crucial role in the expression of $\expect[x]{\zeta_k}$. Hence, by characterising this ratio, we can characterise $\expect[x]{\zeta_k}$. We note from~\eqref{eq:transition_prob} that

\begin{equation}
    \frac{p_{i,i-1}}{p_{i,i+1}}=f_\alpha(i/n),
\end{equation}
where $f_\alpha:[0,1]\to \mb{R}_+$, for each $\alpha \in (0,1)$, is as defined in Theorem~\ref{thm:complete}.
In the lemma below, we obtain a list some important properties of the function $f_\alpha$.

% To prove the theorem above we use a technique developed in~\cite{arpan_JSP_20} where suitably defined branching processes are used to obtain bounds on the number of visits of the chain $\bar Y^n$ to a given state $k \in \cbrac{0,1,\ldots,n}$. Using these bounds and the properties of the function $f_\alpha$ defined in Theorem~\ref{thm:complete} we 

\begin{lemma}
\label{lem:ratio}
For $\alpha \in (0,1)$, define $f_\alpha:[0,1]\to \mb{R}_+$

\begin{equation}
    \label{eq:def_f}
    f_\alpha(x)\define\frac{(1-\alpha)x(1-x)}{\alpha+(1-\alpha)x^2}
\end{equation}
Then $f_\alpha$ satisfies the following properties
\begin{enumerate}
    
    \item For all $\alpha \in (0,1)$, $f_\alpha$ is strictly increasing in $[0,x_\alpha)$, strictly decreasing in $(x_\alpha,  1]$, and, in the domain $[0,1]$, attains its maximum value at $x_\alpha$, where 
\begin{equation}
	x_\alpha\define \sqrt{\brac{\frac{\alpha}{1-\alpha}}^2+\frac{\alpha}{1-\alpha}}-\frac{\alpha}{1-\alpha} \in [0,1).
	\label{eq:x_alpha}
\end{equation}

    \item Let $r_\alpha\define f_\alpha(x_\alpha)=\max_{x \in [0,1]}f_\alpha(x)$. Then, for $\alpha \in (1/9,1)$, we have $r_\alpha < 1$, and, for $\alpha \in (0,1/9)$, we have $r_\alpha > 1$.
    
    \item For $\alpha \in (0,1/9]$, we have $f_\alpha(x) \geq 1$ iff $x \in [\underline{x}_{\alpha},\bar{x}_\alpha]$, where
    
    \begin{align}
        \underline{x}_\alpha &= \frac{1}{4}\brac{1-\sqrt{1-\frac{8\alpha}{1-\alpha}}}\label{eq:x_alpha_ubar}\\
        \bar x_\alpha &=\frac{1}{4}\brac{1+\sqrt{1-\frac{8\alpha}{1-\alpha}}}\label{eq:x_alpha_bar}
    \end{align}
    Furthermore, $f_\alpha(\underline{x}_\alpha)=f_\alpha(\bar{x}_\alpha)=1$ and $x_\alpha \in [\underline{x}_{\alpha},\bar{x}_\alpha]$.
    
    \item For $\alpha \in (0,1/9)$, define $g_\alpha(x)\define\int_{\underline{x}_\alpha}^x \log(f_\alpha(x))dx$.
    Then $g_\alpha$ has a unique root $p_c(\alpha) \in (\bar x_\alpha,1)$. Furthermore, $g_\alpha(p) > 0$ if $p \in [\bar x_\alpha, p_c(\alpha))$ and $g_\alpha(p) < 0$ if $p \in (p_c(\alpha),1)$

    % \item For $\alpha \in (0,1/9]$, we clearly have $\int_{\underline{x}_\alpha}^{\bar x_\alpha} \log(f(x))dx \geq 0$ since $f(x) \geq 1$. It is also easy to see that $\int_{\underline{x}_\alpha}^{1} \log(f(x))dx < 0$ (by directly evaluating the integral). Since $\log(f(x)) \leq 0$ for $x\in [\bar x_\alpha,1]$, there exists a unique $p_c \in [\bar x_\alpha,1)$ such that $\int_{\underline{x}_\alpha}^{p_c} \log(f(x))dx=0$.  
    
\end{enumerate}
\end{lemma}

\begin{proof}
{Taking the derivative of~\eqref{eq:def_f} with respect $x$ we obtain
\begin{equation*}
  f_\alpha'(x)=\frac{-(1-\alpha)^2 x^2-2\alpha(1-\alpha)x+\alpha(1-\alpha)}{\brac{\alpha+(1-\alpha)x^2}^2}.
\end{equation*}
We note that the numerator of the above expression is zero when $x=x_\alpha$, positive when $x \in [0, x_\alpha)$ and negative when $x \in (x_\alpha, 1]$, where $x_\alpha$ is as defined in~\eqref{eq:x_alpha}. This proves the first
statement of the lemma.}

{Next, we note from~\eqref{eq:def_f} that the condition $f_\alpha(x) \lesseqgtr 1$ is equivalent to $2x^2-x+\frac{\alpha}{1-\alpha} \gtreqless 0$.
Hence, $f_\alpha(x) < 1$ for all $x \in [0,1]$
if and only if $1-8\frac{\alpha}{1-\alpha} <0$ or equivalently iff $\alpha > 1/9$. 
%Hence, if $\alpha \in (1/9,1)$, we have $r_\alpha=\max_{x \in [0,1]}f(x) < 1$. 
Furthermore,  for $\alpha \in (0,1/9]$, we have $2x^2-x+\frac{\alpha}{1-\alpha}=(x-\underline{x}_\alpha)(x-\bar x_\alpha)$ where $\underline{x}_\alpha$ and $\bar{x}_\alpha$ are as defined in~\eqref{eq:x_alpha_ubar} and~\eqref{eq:x_alpha_bar}, respectively, and for $\alpha < 1/9$ we have $\underline{x}_\alpha < \bar{x}_\alpha$ and $x_\alpha \in (\underline{x}_\alpha,\bar{x}_\alpha)$.
Combining the above facts we have the second and the third statements of the lemma. 
%Furthermore, if $\alpha < 1/9$, then $\underline{x}_\alpha < x_\alpha < \bar x_\alpha$ and for all $x \in (\underline{x}_\alpha,\bar{x}_\alpha)$ we have $f_\alpha(x) > 1$. The later implies that $r_\alpha=f(x_\alpha) > 1$ if $\alpha \in (0,1/9)$.
}

% The proofs of the first three statements are elementary.
% % \arpan{Add some details here, e.g., expression of the derivative etc.}
% Here, we only prove the last statement of the lemma.

We now turn to the last statement of the lemma. Note that for $\alpha \in (0,1/9)$, we have $$g_\alpha(\bar x_\alpha)=\int_{\underline{x}_\alpha}^{\bar x_\alpha} \log(f_\alpha(x))dx \geq \log(f_\alpha(x_\alpha)) =\log(r_\alpha) > 0,$$ 
where the last inequality follows from the fact that $r_\alpha=f_\alpha(x_\alpha) > 1$ for $\alpha \in (0,1/9)$ as established in the previous paragraph.
% since by the third statement of the lemma  $x_\alpha \in [\underline{x}_\alpha,\bar x_\alpha]$ for $\alpha \in (0,1/9]$ and by the second statement of the lemma $r_\alpha > 1$ for $\alpha \in (0,1/9)$. 
Using~\eqref{eq:def_f}, we can compute $g_\alpha$ in closed form. This is given as follows

\begin{align}
     g_\alpha(x) &=
    %  x\log\brac{\frac{f(x)}{1-x}}-\underline{x}_\alpha \log(f(\underline{x}_\alpha)) - (1-x)\log(1-x)+\log(1-\underline{x}_\alpha) \nonumber\\
    % & -2\sqrt{\frac{\alpha}{1-\alpha}}\brac{\arctan\brac{\sqrt{\frac{1-\alpha}{\alpha}}x}-\arctan\brac{\sqrt{\frac{1-\alpha}{\alpha}}\underline{x}_\alpha}} \\
     x\log\brac{\frac{(1-\alpha)x}{\alpha+(1-\alpha)x^2}}- (1-x)\log(1-x)+\log(1-\underline{x}_\alpha) \nonumber\\
    & -2\sqrt{\frac{\alpha}{1-\alpha}}\brac{\arctan\brac{\sqrt{\frac{1-\alpha}{\alpha}}x}-\arctan\brac{\sqrt{\frac{1-\alpha}{\alpha}}\underline{x}_\alpha}}
\end{align}
From the above it follows that $\lim_{x \to 1^{-}} g(x) < 0$
since $\underline{x}_\alpha < 1$ and $\arctan(\cdot)$ is an increasing function.
Furthermore, we have $g'_\alpha(x)=\log(f_\alpha(x)) < 0$ for $x\in (\bar x_\alpha,1)$. Hence, there must exist a unique root $p_c(\alpha)$ of $g_\alpha$
in $(\bar x_\alpha, 1)$ and $g_\alpha(p)$ must be strictly negative for $p \in (p_c(\alpha),1)$ and strictly positive for $p \in [\bar{x}_\alpha,p_c(\alpha))$.
\end{proof}

We are now in a position to complete the proof of Theorem~\ref{thm:complete}. We use Lemma~\ref{lem:zeta_k} and the properties of $f_\alpha$ proved in Lemma~\ref{lem:ratio} to obtain upper bounds of $\expect[x]{\zeta_k}$. These upper bounds, in turn, provide upper bounds on $\expect[x]{Z_k}$ using~\eqref{eq:z_k}. Finally, we use Lemma~\ref{lem:master} and the upper bounds on $\expect[x]{Z_k}$ to obtain an upper bound on the mean consensus time $\bar T_n(p)$.
The complete proof is given below.

{\em Proof of Theorem~\ref{thm:complete}}: From Lemma~\ref{lem:master} it follows that 
% Now we show that $\bar T_n(p)=O(n \log n)$ under the conditions stated in the first
% two parts of the theorem. To do so, we first 
% note the following from~\eqref{eq:consensus}

\begin{align}
    \bar T_n(p) &=\sum_{k=0}^{n-1} \frac{\expect[\ceil{pn}]{Z_k}}{(1-p_{k,k})} \nonumber\\
    &= \sum_{k=0}^{n-1} \frac{\expect[\ceil{pn}]{Z_k}}{(1-\frac{k}{n})(\alpha+(1-\alpha)\frac{k}{n})} \nonumber\\
    &=n\sum_{k=0}^{n-1} \expect[\ceil{pn}]{Z_k}\brac{\frac{1}{n-k}+\frac{1}{\frac{\alpha}{1-\alpha}n+k}} \label{eq:tnp}
\end{align}
{Hence, to prove $T_n(p)=O(n \log n)$, it is sufficient to establish that $\expect[\ceil{p n}]{Z_k} \leq C$ for all $k \in \{0,1,\ldots,n\}$, where $C >0$ is a constant independent of $k$. Indeed, if $\expect[\ceil{p n}]{Z_k} \leq C$ for each state $k \in \{0,1,\ldots,n\}$, then for $n \geq \frac{1-\alpha}{\alpha}$
we have}

\begin{align}
    \bar T_n(p) & \leq n \textcolor{blue}{C}\sum_{k=0}^{n-1} \brac{\frac{1}{n-k}+\frac{1}{k+1}} \nonumber \\
    &=2 n \textcolor{blue}{C}\sum_{k=1}^{n-1}\frac{1}{k} \nonumber\\
    &=O(n\log n)
\end{align}
Hence, to prove the first two statements of the theorem, it is sufficient to show that, under the conditions stated in these statements, $\expect[\ceil{pn}]{Z_k}$ is uniformly bounded by some constant for all $k \in \cbrac{0,1,\ldots,n-1}$. This what we prove next.

Using the first two properties of $f_\alpha$ proved in
Lemma~\ref{lem:ratio} and~\eqref{eq:solved_expect_zeta},
we see that for
$\alpha \in (1/9,1)$ we have

\begin{equation}
\label{eq:ineq_expect_zeta}
\expect[x]{\zeta_{k}}\leq \begin{cases}
			\sum_{t=k}^{n-1} r_\alpha^{t-k+1} < \frac{r_\alpha}{1-r_\alpha}, &\quad \text{for } x \leq k \leq n-1\\
			\expect[x]{\zeta_{x}} < \frac{r_\alpha}{1-r_\alpha}, &\quad \text{for } k < x,
			\end{cases}
\end{equation}
where $r_\alpha < 1$ is as defined in Lemma~\ref{lem:ratio}.
Hence, from~\eqref{eq:z_k} we have that for all $\alpha \in (1/9,1)$ and all $k \in \cbrac{0,1,\ldots,n-1}$, $$\expect[x]{Z_k}\leq  \frac{1+r_\alpha}{1-r_\alpha}.$$
This establishes the first statement of the theorem.

We now prove the second statement of the theorem. For $\alpha \in (0,1/9)$, the third statement of Lemma~\ref{lem:ratio} implies $p_c(\alpha) > \bar x_\alpha$.
Furthermore, by the first and third statements of Lemma~\ref{lem:ratio}, it follows that $f_\alpha$ is strictly decreasing in the range $[\bar x_\alpha,1)$. Hence, for $p \geq p_c(\alpha) > \bar x_\alpha$ we have
$f_\alpha(p) < f_\alpha(\bar x_\alpha) =1$. 
Furthermore, for all $k \geq x \geq pn$ we have
$$\frac{p_{k,k-1}}{p_{k,k+1}}=f_\alpha(k/n)\leq f_\alpha(p) < 1.$$
Let $r_p\define f_\alpha(p) < 1$. Then, from~\eqref{eq:solved_expect_zeta}, we have 
that for $k \geq x$
$$\expect[x]{\zeta_k}\leq \sum_{t=k}^{n-1} r_p^{t-k+1} < \frac{r_p}{1-r_p}.$$
Hence,
from \eqref{eq:z_k} we have
that $$\expect[x]{Z_k} \leq \frac{1+r_p}{1-r_p}$$
for $k \geq x$. 
Moreover, for $k<  n\bar{x}_\alpha \leq x= \ceil{pn}$ we have

\begin{align*}
    \prod_{i=k}^{x-1} \frac{p_{i,i-1}}{p_{i,i+1}}
    &=\prod_{i=k}^{\ceil{pn}-1} f_\alpha(i/n)\\
    & \overset{(a)}{\leq} \prod_{i=\ceil{n\underline{x}_\alpha}}^{\ceil{pn}-1} f(i/n)\\
    &=\exp\brac{\sum_{i=\ceil{n\underline{x}_\alpha}}^{\ceil{pn}-1}\log\brac{f_\alpha(i/n)}}\\
    &{\overset{(b)}{=}\exp\brac{n\brac{ \int_{\underline{x}_\alpha}^{p} \log (f_\alpha(x))dx +O(1/n)}}}\\
    & \overset{(c)}{=}\exp\brac{n {g_\alpha(p)} +O(1)}\\
    &\overset{(d)}{=}O(1),
\end{align*}
where (a) follows from the facts that
$f_\alpha(i/n) \leq 1$ for $i \leq n\underline{x}_\alpha$ and $f_\alpha(i/n) \geq 1$ for $n\underline{x}_\alpha \leq i \leq n\bar{x}_\alpha$; (b) follows from the fact that the Riemannian sum $(1/n)\sum_{i=\ceil{n\underline{x}_\alpha}}^{\ceil{pn}-1}\log\brac{f_\alpha(i/n)}$ converges
to the integral $\int_{\underline{x}_\alpha}^{p} \log (f_\alpha(x))dx$ as $n \to \infty$ with an error of $O(1/n)$; (c) follows from the definition of $g_\alpha$ in Lemma~\ref{lem:ratio}; (d) follows from the fact
that for $p \geq p_c(\alpha) > \bar x_\alpha$ we have
$g_\alpha(p) \leq 0$ by the last statement of Lemma~\ref{lem:ratio}. For $x > k \geq n \bar x_\alpha$
we have 

\begin{align*}
    \prod_{i=k}^{x-1} \frac{p_{i,i-1}}{p_{i,i+1}}
    &=\prod_{i=k}^{\ceil{pn}-1} f(i/n)
    \leq 1,
\end{align*}
because for $i \geq n \bar{x}_\alpha$
we have $f_\alpha(i/n) \leq 1$.
Hence, for all $k < x$, it follows from 
\eqref{eq:solved_expect_zeta} that
$\expect[x]{\zeta_k}=O(1) \expect[x]{\zeta_x}=O(1)$.
Finally, from \eqref{eq:z_k} we obtain
that $\expect[x]{Z_k} =O(1)$ for all $k < x$.
This establishes the second statement of the theorem.

We now turn to the third statement of the theorem. Note from~\eqref{eq:tnp} that to prove this statement, it suffices to show that $\expect[\ceil{pn}]{Z_k}=\Omega(\exp(\Theta(n)))$ for $k =\ceil{n\underline{x}_\alpha}$ when $\alpha \in (0,1/9)$ and $p \in (0, p_c(\alpha))$. First, note that
for $\alpha \in (0,1/9)$, $p \in (\underline{x}_\alpha,p_c(\alpha))$,
and $k=\ceil{n\underline{x}_\alpha}$ we have 
\begin{align*}
\prod_{i=\ceil{n\underline{x}_\alpha}}^{x-1} \frac{p_{i,i-1}}{p_{i,i+1}}
&=\prod_{i=\ceil{n\underline{x}_\alpha}}^{\ceil{np}-1} f_\alpha(i/n)\\
&=\exp\brac{n \int_{\underline{x}_\alpha}^{p}\log(f_\alpha(x))dx +O(1)}\\
&=\Omega(\exp(\Theta(n))),
\end{align*}
where the last line follows from the fact that 
$\int_{\underline{x}_\alpha}^{p}\log(f_\alpha(x))dx=g_\alpha(p) > 0$ for $p\in (\underline{x}_\alpha,p_c(\alpha))$.
Hence, $\expect[x]{\zeta_{\ceil{\underline{x}_\alpha n}}}=\Omega(\exp(\Theta(n)))\expect[x]{\zeta_x}
=f_\alpha(p)\Omega(\exp(\Theta(n)))$ since
$\expect[x]{\zeta_x}>f_\alpha(x/n)=f_\alpha(p)+o(1)$.
Hence, $\expect[x]{Z_{\ceil{n\underline{x}_\alpha}}}=\Omega(\exp(\Theta(n)))$.
For $p \in [0,\underline{x}_\alpha]$ we have
\begin{align*}
\expect[x]{\zeta_{\ceil{\underline{x}_\alpha n}}}> \prod_{i=\ceil{n\underline{x}_\alpha}}^{\ceil{n\bar{x}_\alpha}-1} \frac{p_{i,i-1}}{p_{i,i+1}}
&=\prod_{i=\ceil{n\underline{x}_\alpha}}^{\ceil{n\bar{x}_\alpha}-1} f_\alpha(i/n)\\
&=\exp\brac{n \int_{\underline{x}_\alpha}^{\bar{x}_\alpha}\log(f_\alpha(x))dx +O(1)}\\
&=\Omega(\exp(\Theta(n))),
\end{align*}
This proves the third part of the theorem.\qed

\begin{remark}
From the proof above, it is clear that the absorption time of the chain $\bar Y^n$ crucially depends on the ratio
 of the down-transition probability $p_{i,i-1}$ to the up-transition probability $p_{i,i+1}$
at any given state $i$.
If this ratio is  smaller than one at a given state $i$, then the chain has a  tendency to
drift quickly towards the absorbing state $n$.   Similarly, if the ratio is larger than one, then
the chain has a tendency to move towards state $0$ and hence it takes longer to reach the absorbing state $n$.  Since the value of this ratio
depends on the bias parameter $\alpha$ as well as on the initial proportion $p$ of agents with the superior opinion, we observe phase transitions with respect to both the parameters.
\end{remark}

\begin{remark}
\label{rem:mod_chain}
Although we have proved the theorem for the modified chain $\bar{Y}^n$,  the proof can be extended to the original chain $\bar X^n$. This is because
the results of Lemma~\ref{lem:master} and Lemma~\ref{lem:zeta_k} are also applicable to the original chain if the  transition probabilities are replaced by
those of the original chain.
Furthermore, the transition probabilities of the two chains satisfy the following relations $$1-\tilde{p}_{k,k}=1-p_{k,k}+(1-\alpha)\frac{k}{n}\brac{1-\frac{k}{n}}\brac{\frac{n^2}{(n-1)^2}-1}\leq 1-p_{k,k},$$ and $\tilde{p}_{k,k-1}/\tilde{p}_{k,k+1}={p_{k,k-1}}/{p_{k,k+1}}+o(1)$,
where we recall that $\Tilde{p}_{k,k}$ denotes the transition probability from state $k$ to itself for the original chain $\bar X^n$.
%Furthermore,  the result in Lemma~\ref{lem:master} is valid in general for any random walk on $\cbrac{0,1,\ldots,n}$ with absorbing state $n$. Hence, we can apply the result to the original chain $\bar X^n$ and write the expected consensus time $\Tilde{T}_n(p)$ for the original chain as 
% $\Tilde{T}_n(p)=\sum_{k=0}^{n-1} \frac{\expect[x]{\tilde{Z} _k}}{(1-\tilde{p}_{k,k})}$, where $\tilde Z_k$ denotes the number of visits to state $k$
% in the original chain $\bar X^n$. 
We note that the same lower bound as in Lemma~\ref{lem:lower_bound} can be obtained for the chain $\bar X^n$ 
since we have $1-\tilde{p}_{k,k} \leq (\max(\alpha,1-\alpha)+o(1))(1-k/n)(1+k/n)$. 
Similarly, the upper and lower bounds on the expected number of visits to each state derived for the modified chain also hold for the original chain for large enough $n$.  Combining the above, we obtain the same asymptotic bounds for the original chain $\bar X^n$ as we have for the modified chain $\bar Y^n$.
%    Since $\tilde p_{k,k}=p_{k,k}+o(1/n)$, we have $\Tilde{T}_n(p)=\Bar{T}_n(p)+o(1)$. This validates the 
\end{remark}

\section{Numerical Results}
\label{sec:numerics}

In this section, we present numerical results for the model presented in this paper. We first present results to support our theoretical findings on complete graphs. We then present simulation results for other classes of graphs.
Error-bars in the plots represent $95\%$ confidence intervals. Also, to understand the growth rates better, we overlay the plots obtained from simulations on theoretical growth rates obtained by plotting appropriately scaled functions.

\subsection{Complete graphs}

\begin{figure}[htb]
     \centering
     \begin{subfigure}[b]{0.45\textwidth}
         \centering
         \includegraphics[width=\textwidth]{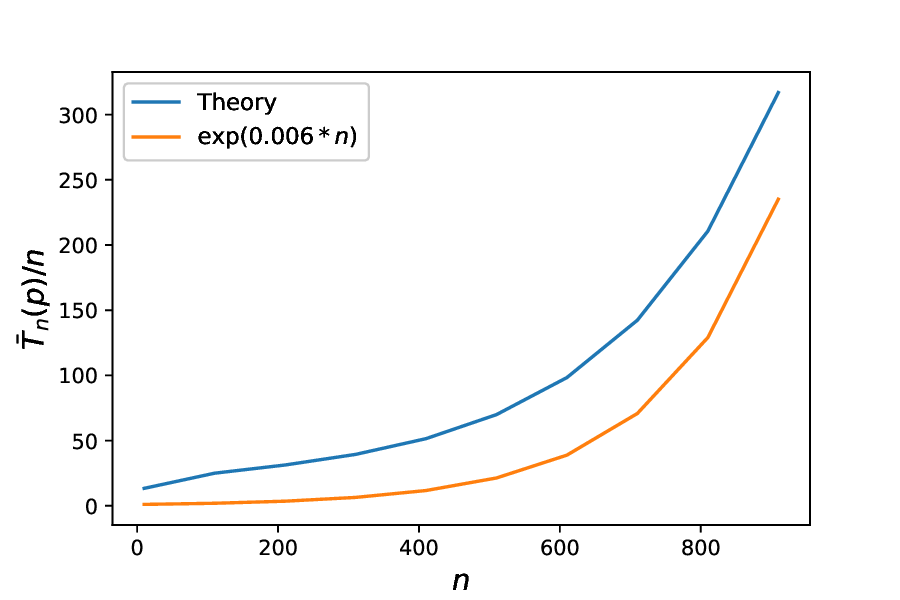}
         \caption{$\alpha=0.1 < 1/9, p=0.4 < p_c(\alpha)=0.431$}
         \label{fig:complete_belowcutoff}
     \end{subfigure}
     \hfill
     \begin{subfigure}[b]{0.45\textwidth}
         \centering
         \includegraphics[width=\textwidth]{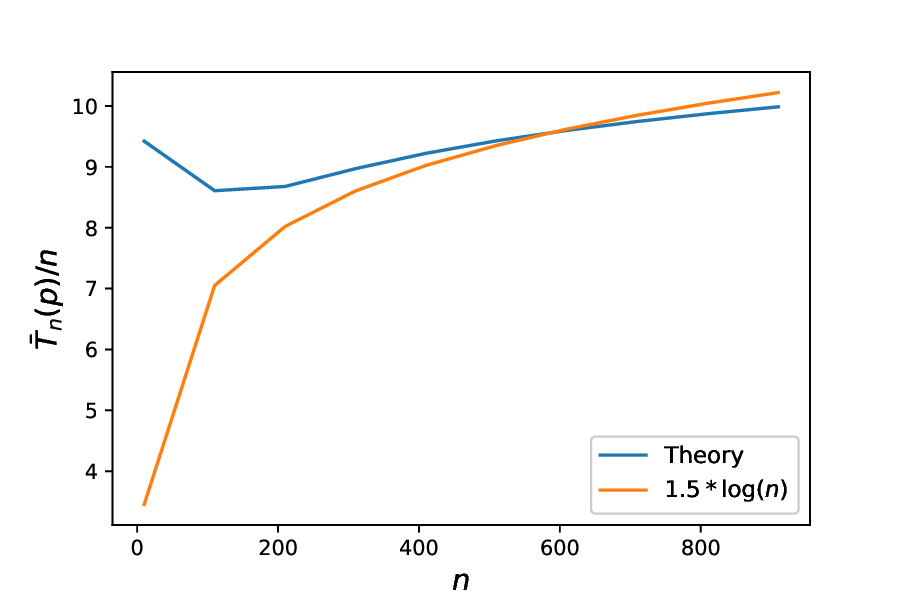}
         \caption{$\alpha=0.1 < 1/9, p=0.5 > p_c(\alpha)=0.431$}
         \label{fig:complete_abovecutoff}
     \end{subfigure}
     \caption{Mean consensus time per node $\bar{T}_n(p)/n$ as a function of the network size for $\alpha=0.1 < 1/9$ for complete graphs}
\end{figure}

{For complete graphs, the mean consensus time $\bar T_n(p)$ can be computed numerically using the first step analysis of the Markov chain $\bar Y^n$. 
This method is more exact and computationally less expensive than simulating the chain a large number of times to obtain the average absorption time.  Hence, we adopt this method for complete graphs. 
With a slight abuse of notation let $\bar T_n(k)$ denote the average absorption time of the chain starting from state $k \in \cbrac{0,1,\ldots,n}$. Then, from first step analysis
we have for each $k \in \cbrac{0,1,\ldots,n}$
\begin{equation}
\bar T_n(k)=1+p_{k,k-1}\bar{T}_n(k-1)+p_{k,k+1}\bar{T}_n(k+1)+p_{k,k}\bar{T}_n(k),
\end{equation}
where $p_{i,j}$ is as defined in~\eqref{eq:transition_prob}.
Simplifying the above and using the boundary condition $\bar T_n(n)=0$ we obtain}

\begin{equation}
    \bar{T}_n(p):= \bar T_n(\ceil{np})=\sum_{k=\ceil{np}}^{n-1}S_n(k),
\end{equation}
where $S_n(k)$ satisfies the following recursion:

\begin{equation}
    S_n(k)=\frac{1}{p_{k,k+1}}+\frac{p_{k,k-1}}{p_{k,k+1}}S_n(k-1),
\end{equation}
with $S_n(0)=1/\alpha$. 
We find the expected consensus time numerically by solving the above recursion for different values of $n, \alpha,$ and $p$. First, we choose $\alpha=0.1 < 1/9$. For this value of $\alpha$, we can numerically compute the value of $p_c(\alpha)$ by solving~\eqref{eq:p_cutoff}. This value turns out to be $0.431$ (accurate to the third decimal place). In Figures~\ref{fig:complete_belowcutoff} and~\ref{fig:complete_abovecutoff} we plot the normalised (by the network size) average consensus time as a function of the network size for $p=0.4 < p_c(\alpha)=0.431$ and $p=0.5 > p_c(\alpha)=0.431$, respectively.
As expected from Theorem~\ref{thm:complete} we observe that the mean consensus time grows exponentially for $p < p_c(\alpha)$ and as $\Theta(n \log n)$ for $p > p_c(\alpha)$.

\begin{figure}[htb]
     \centering
     \begin{subfigure}[b]{0.45\textwidth}
         \centering
         \includegraphics[width=\textwidth]{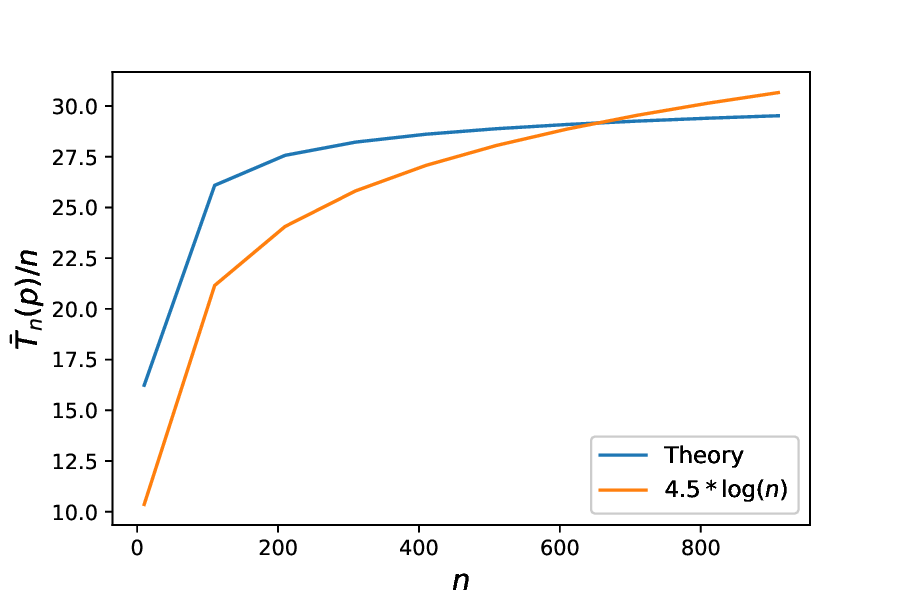}
         \caption{$\alpha=0.125 > 1/9, p=0$}
         \label{fig:complete_p0}
     \end{subfigure}
     \hfill
     \begin{subfigure}[b]{0.45\textwidth}
         \centering
         \includegraphics[width=\textwidth]{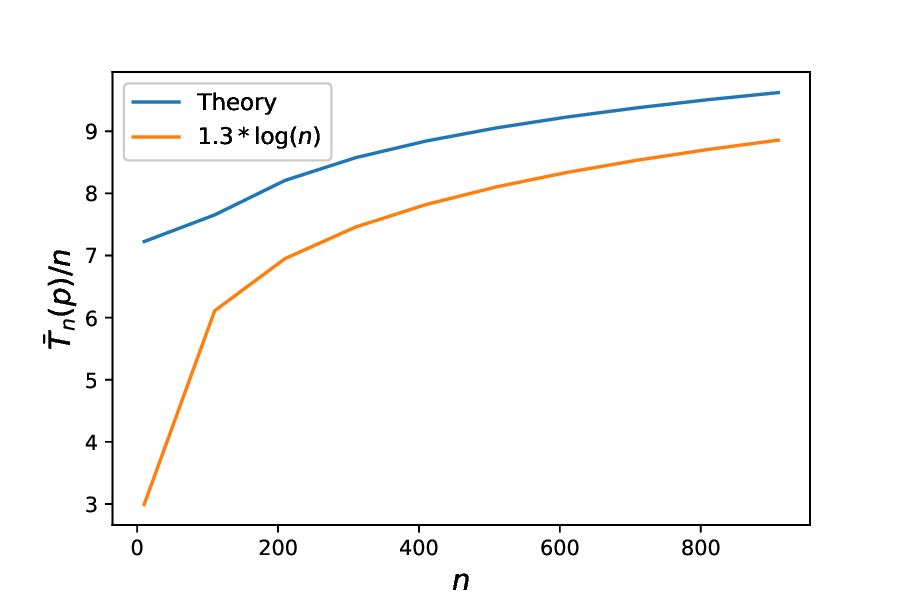}
         \caption{$\alpha=0.125 > 1/9, p=0.5$}
         \label{fig:complete_p05}
     \end{subfigure}
     \caption{Mean consensus time per node $\bar{T}_n(p)/n$ as a function of the network size for $\alpha=0.125 > 1/9$ for complete graphs}
\end{figure}

Next, we choose $\alpha =0.125 > 1/9$. For this choice of $\alpha$, we expect (from Theorem~\ref{thm:complete}) the mean consensus time to grow as $\Theta(n \log n)$ for all values of $p$. This is verified in Figures~\ref{fig:complete_p0} and~\ref{fig:complete_p05}, where we choose $p=0$ and $p=0.5$, respectively. We observe that in both cases $\bar{T}_n(p)=\Theta(n \log n)$.

\subsection{Random $d$-regular graphs with $d=\Theta(\log n)$}

\begin{figure}[htb]
     \centering
     \begin{subfigure}[b]{0.45\textwidth}
         \centering
         \includegraphics[width=\textwidth]{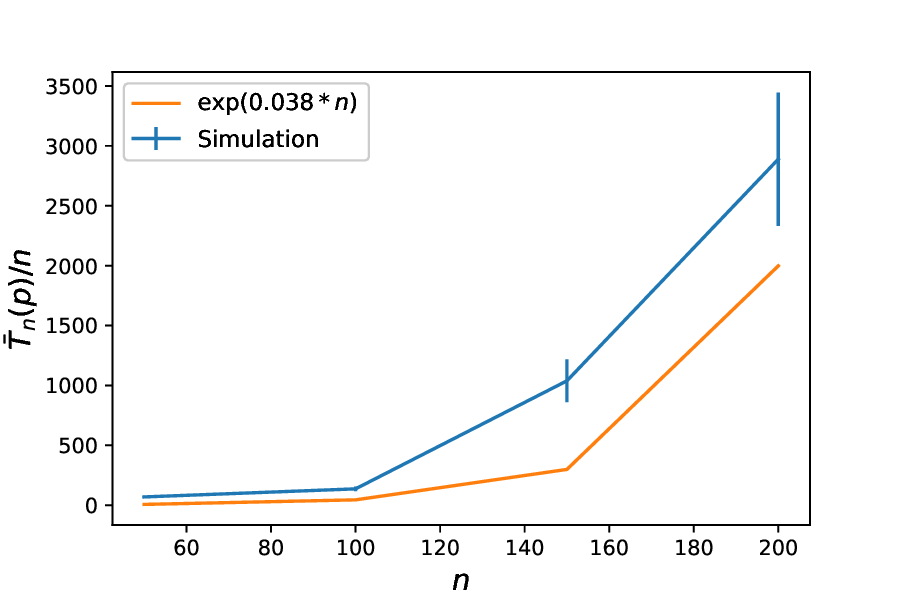}
         \caption{$\alpha=0.05, p=0.05$}
         \label{fig:regular_log_alpha005_p005}
     \end{subfigure}
     \hfill
     \begin{subfigure}[b]{0.45\textwidth}
         \centering
         \includegraphics[width=\textwidth]{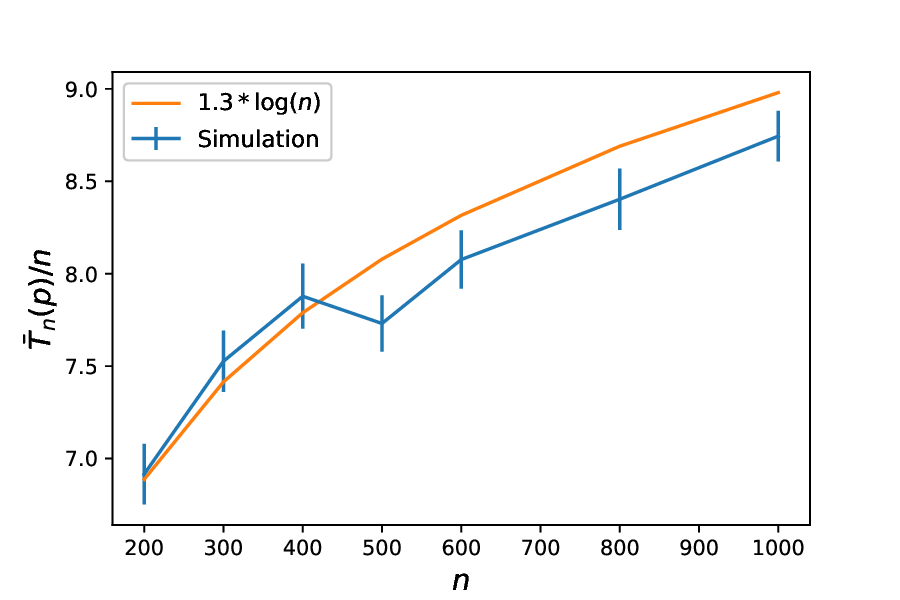}
         \caption{$\alpha=0.05, p=0.8$}
         \label{fig:regular_log_alpha005_p08}
     \end{subfigure}
     \caption{Mean consensus time per node $\bar{T}_n(p)/n$ as a function of the network size for $\alpha=0.05$ for random $d$-regular graphs with $d=\ceil{\log n}$.}
\end{figure}

We now present simulation results for random $d$-regular graphs where the degree $d$ for each node is chosen to be $d=\ceil{\log n}$. We use the $\texttt{networkx}$ package to generate random $d$-regular graphs. For each randomly generated graph, we run the protocol until the network reaches consensus. The above procedure is repeated $500$ times and mean consensus time is computed over these $500$ runs. We further repeat this for different values of $n,\alpha$, and $p$. In Figures~\ref{fig:regular_log_alpha005_p005} and~\ref{fig:regular_log_alpha005_p08} we fix $\alpha$ to be $0.05$ and plot the normalised mean consensus time as a function of the network size for $p=0.05$ and $p=0.8$, respectively.  Similar to complete graphs, we observe that the mean consensus time grows exponentially with $n$ when both $\alpha$ and $p$ are low. In contrast, when the initial proportion of agents having the superior opinion is sufficiently large, the mean consensus time grows as $\Theta(n\log n)$ even when the bias parameter $\alpha$ is small. 

\begin{figure}[htb]
     \centering
     \begin{subfigure}[b]{0.45\textwidth}
         \centering
         \includegraphics[width=\textwidth]{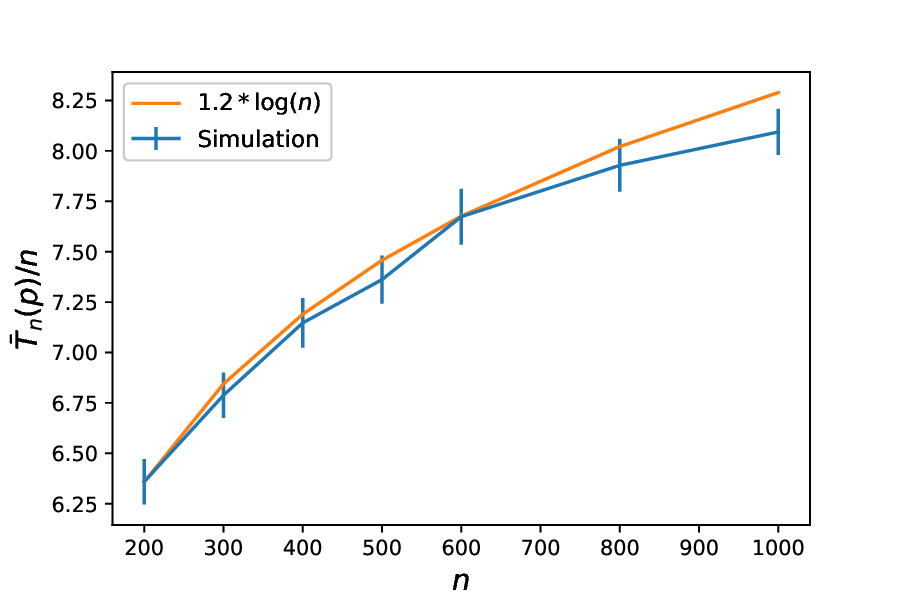}
         \caption{$\alpha=0.8, p=0.05$}
         \label{fig:regular_log_alpha08_p005}
     \end{subfigure}
     \hfill
     \begin{subfigure}[b]{0.45\textwidth}
         \centering
         \includegraphics[width=\textwidth]{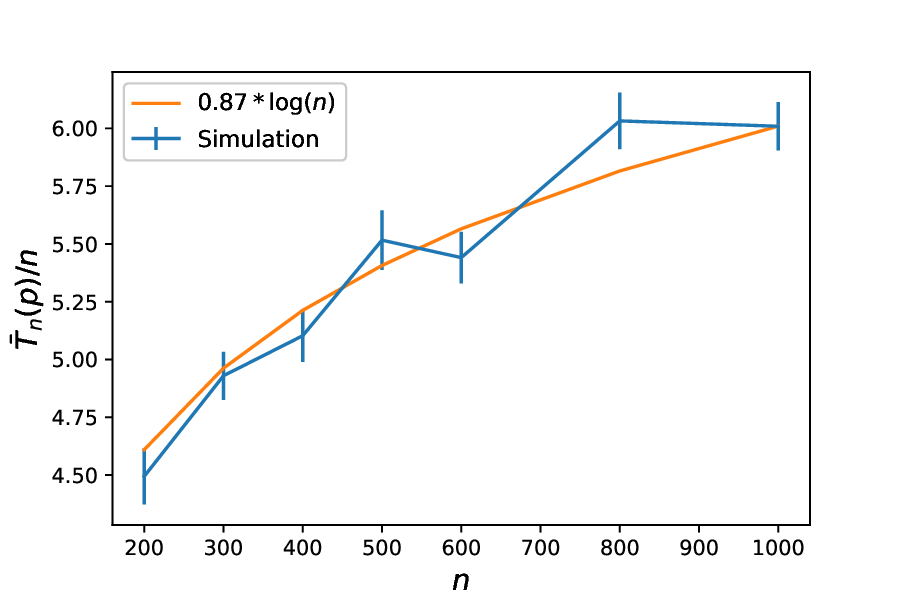}
         \caption{$\alpha=0.8, p=0.8$}
         \label{fig:regular_log_alpha08_p08}
     \end{subfigure}
     \caption{Mean consensus time per node $\bar{T}_n(p)/n$ as a function of the network size for $\alpha=0.8$ for random $d$-regular graphs with $d=\ceil{\log n}$.}
\end{figure}

In Figures~\ref{fig:regular_log_alpha08_p005} and~\ref{fig:regular_log_alpha08_p08} we choose $\alpha=0.8$ plot the normalised mean consensus time as a function of the network size for $p=0.005$ and $p=0.8$, respectively. We observe that in both cases the mean consensus time grows as $\Theta(n \log n)$. This indicates that if the bias parameter $\alpha$ is sufficiently large then the network reaches consensus to the superior opinion in $\Theta(n \log n)$ time irrespective of the initial proportion of agents having the superior opinion. 

\subsection{Random $d$-regular graphs with $d=O(1)$}

\begin{figure}[htb]
     \centering
     \begin{subfigure}[b]{0.45\textwidth}
         \centering
         \includegraphics[width=\textwidth]{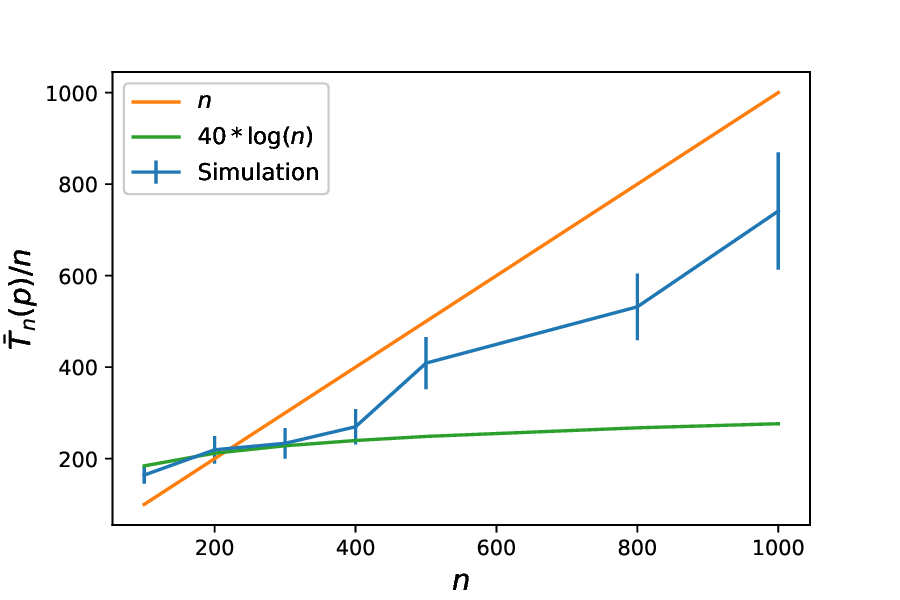}
         \caption{$\alpha=0.05, p=0.05$}
         \label{fig:regular_constant_alpha005_p0}
     \end{subfigure}
     \hfill
     \begin{subfigure}[b]{0.45\textwidth}
         \centering
         \includegraphics[width=\textwidth]{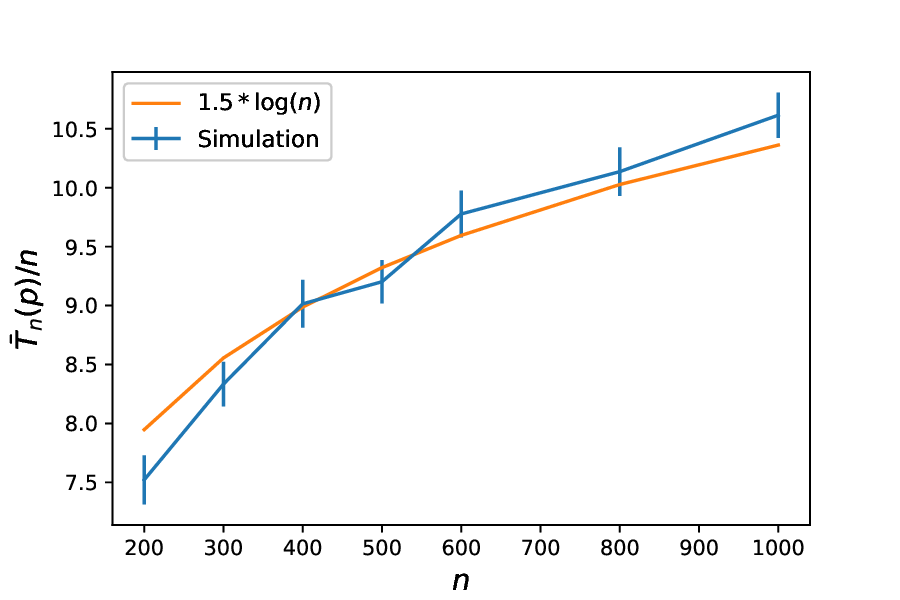}
         \caption{$\alpha=0.05, p=0.8$}
         \label{fig:regular_constant_alpha005_p08}
     \end{subfigure}
     \caption{Mean consensus time per node $\bar{T}_n(p)/n$ as a function of the network size for $\alpha=0.05$ for random $d$-regular graphs with $d=5$.}
\end{figure}

We now investigate the biased opinion dynamics on random $d$-regular with constant degree, i.e., $d=O(1)$.
Specifically, we set $d=5$ and simulate the network for different values of $n,\alpha,$ and $p$. As before, we generate a random $5$-regular graph using the \texttt{networkx} package and run the protocol on this randomly generated instance until consensus is reached. The above procedure is repeated multiple times (each with a newly generated random graph) to obtain the mean consensus time within the desired confidence bounds. In Figures~\ref{fig:regular_constant_alpha005_p0} and~\ref{fig:regular_constant_alpha005_p08}, we fix $\alpha=0.05$ and choose $p=0.05$ and $p=0.8$, respectively. We observe that the normalised mean consensus time grows approximately linearly in $n$ (i.e., $\bar T_n(p)= O(n^2)$) when both $\alpha$ and $p$ are small. This is unlike the previous results for dense graphs, where, for small values of $\alpha$ and $p$, the mean consensus time is exponential in $n$.
This reduction in absorption time can be intuitively explained by the fact that a smaller neighbourhood size reduces the chance of an agent updating to the worse opinion when a fixed proportion of all the agents have the worse opinion. 
Based on this intuition we conjecture that for small values of $p$ and $\alpha$
the mean consensus time grows polynomially in $n$ when the neighbourhood size $d$ is constant higher than $2$.
%Thus, when the neighbourhood sizes are small, consensus on the superior opinion is reached quicker. Intuitively, this means that smaller neighbourhood sizes reduce social pressure and help achieve consensus faster when both the bias parameter and the initial proportion of agents with the superior opinion are low. 
{The case of cycles (where $d=2$) has already been considered in~\cite{biased_info_sciences_2022}. 
But in the case of cycles  there is no phase transition; indeed,  it has been shown in~\cite{biased_info_sciences_2022}
that for cycles consensus is achieved on the superior opinion in $O(n\log n)$ time for all values of $\alpha$ even when $p=0$.}

\begin{figure}[htb]
     \centering
     \begin{subfigure}[b]{0.45\textwidth}
         \centering
         \includegraphics[width=\textwidth]{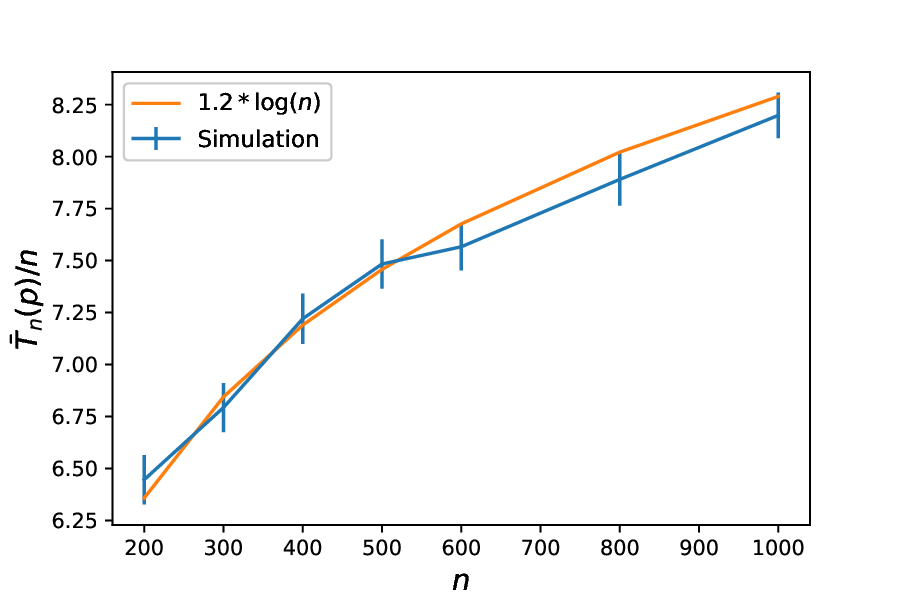}
         \caption{$\alpha=0.8, p=0.05$}
         \label{fig:regular_constant_alpha08_p005}
     \end{subfigure}
     \hfill
     \begin{subfigure}[b]{0.45\textwidth}
         \centering
         \includegraphics[width=\textwidth]{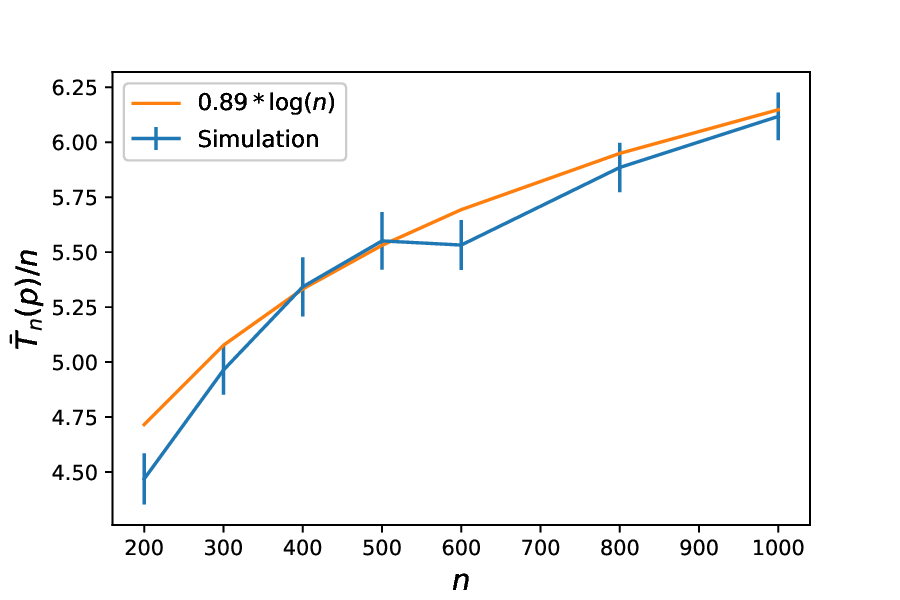}
         \caption{$\alpha=0.8, p=0.8$}
         \label{fig:regular_constant_alpha08_p08}
     \end{subfigure}
     \caption{Mean consensus time per node $\bar{T}_n(p)/n$ as a function of the network size for $\alpha=0.8$ for random $d$-regular graphs with $d=5$.}
\end{figure}

In Figures~\ref{fig:regular_constant_alpha08_p005} and~\ref{fig:regular_constant_alpha08_p08}, we plot the mean consensus time as a function of the network size for $p=0.05$ and $p=0.8$, respectively, keeping $\alpha=0.8$. In both cases, we observe that the mean consensus time grows as $\Theta(n \log n)$. 
Thus, in these regimes, the behaviour of random $d$-regular graphs with $d=O(1)$ is similar to that of random $d$-regular graphs with $d=\Theta(\log n)$.
% The above results indicate that the phase transition we observe for dense graphs (where the degree of each node is $d=\Omega(\log n)$) is absent when the graph is sparse, i.e., when $d=O(1)$. This result was proved for the special case of cycles (where $d=2$) in~\cite{biased_info_sciences_2022}. However, a general proof remains an open problem.

\subsection{\ER graphs with edge probability $\log n/n$}

\begin{figure}[htb]
     \centering
     \begin{subfigure}[b]{0.45\textwidth}
         \centering
         \includegraphics[width=\textwidth]{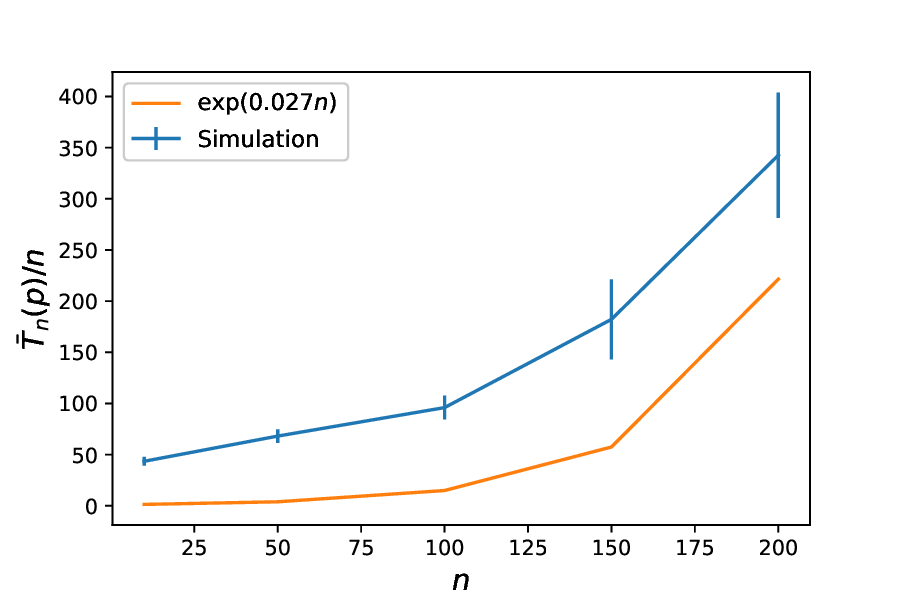}
         \caption{$\alpha=0.05, p=0.05$}
         \label{fig:ER_log_alpha005_p005}
     \end{subfigure}
     \hfill
     \begin{subfigure}[b]{0.45\textwidth}
         \centering
         \includegraphics[width=\textwidth]{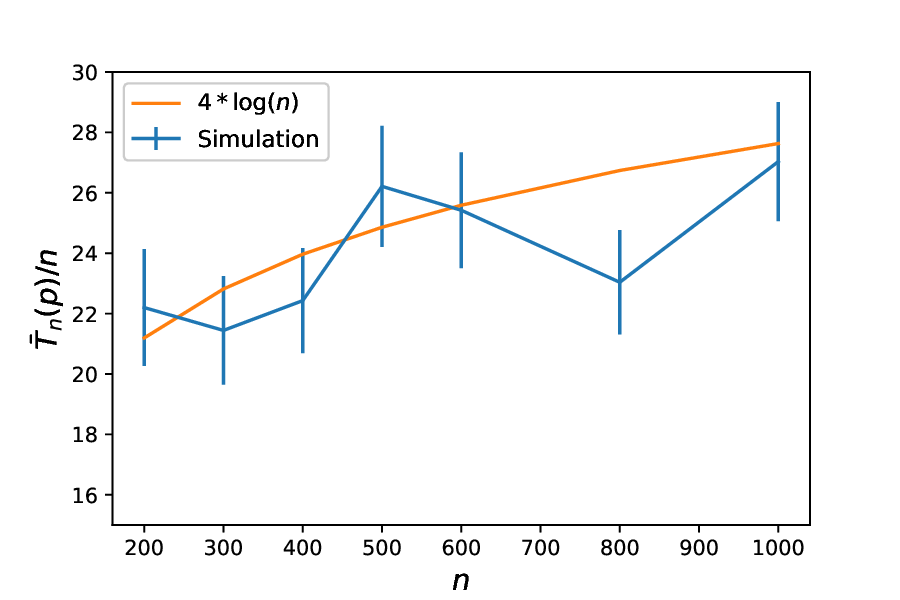}
         \caption{$\alpha=0.05, p=0.8$}
         \label{fig:ER_log_alpha005_p08}
     \end{subfigure}
     \caption{Mean consensus time per node $\bar{T}_n(p)/n$ as a function of the network size for $\alpha=0.05$ for \ER graphs with edge probability $\log n/n$.}
\end{figure}

\begin{figure}[htb]
     \centering
     \begin{subfigure}[b]{0.45\textwidth}
         \centering
         \includegraphics[width=\textwidth]{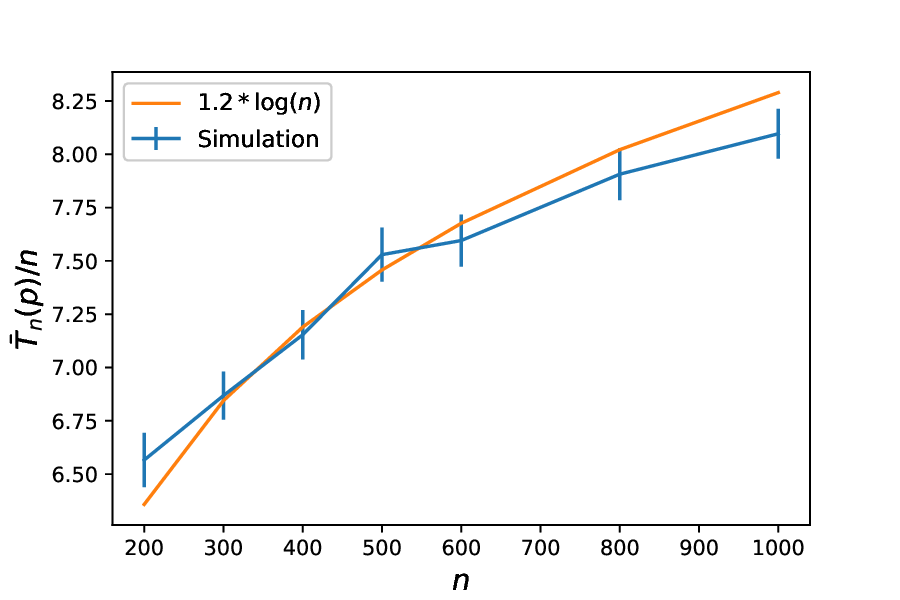}
         \caption{$\alpha=0.8, p=0.05$}
         \label{fig:ER_log_alpha08_p005}
     \end{subfigure}
     \hfill
     \begin{subfigure}[b]{0.45\textwidth}
         \centering
         \includegraphics[width=\textwidth]{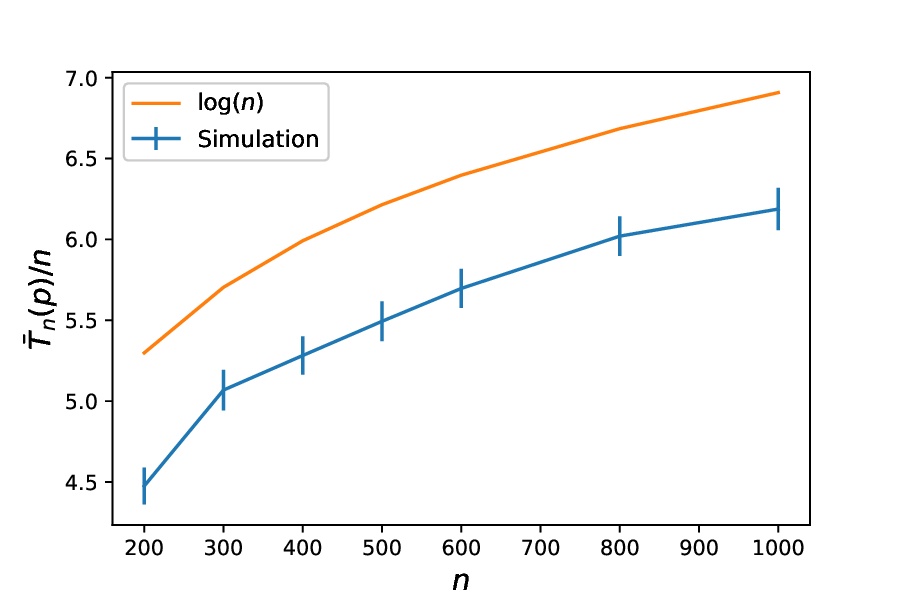}
         \caption{$\alpha=0.8, p=0.8$}
         \label{fig:ER_log_alpha08_p08}
     \end{subfigure}
     \caption{Mean consensus time per node $\bar{T}_n(p)/n$ as a function of the network size for $\alpha=0.8$ for \ER graphs with edge probability $\log n/n$.}
\end{figure}

Finally, we present results for \ER graphs where {the edge} probability is fixed at the connectivity threshold $\log n/n$. The experiments are designed in the same way as described for other classes of graphs. The normalised mean consensus time as a function of the network size $n$ is plotted for different values of $\alpha$ and $p$: in Figures~\ref{fig:ER_log_alpha005_p005} and~\ref{fig:ER_log_alpha005_p08} for $\alpha=0.05$ and in Figures~\ref{fig:ER_log_alpha08_p005} and~\ref{fig:ER_log_alpha08_p08} for $\alpha=0.8$. We observe similar phase transitions as in the case of complete graphs, i.e., the mean consensus time grows as $\Theta(n\log n)$ in all cases except when both $\alpha$ and $p$ are low.

\section{Conclusion and future work}
\label{sec:conclusion}

In this paper we have studied a model of binary opinion dynamics where the agents show a strong form of bias towards one of the opinions, called the superior opinion. We showed that for complete graphs the model exhibits rich phase transitions based on the values of the bias parameter and the initial proportion of agents with the superior opinion. Specifically, we proved that fast consensus can be achieved on the superior opinion irrespective of the initial configuration of the network when bias is sufficiently high.
When bias is low, we show that fast consensus can only be achieved when the initial proportion of agents with the superior opinion is above a certain threshold. If this is not the case, then we show that consensus takes exponentially long time. Through simulations, we observed similar behaviour for several classes of dense graphs where the average degree scales at least logarithmically with the network size. For sparse graphs, where the average degree is constant, we observed that phase transitions do occur but the behaviour below criticality is different to that of dense graphs. Specifically, we observed that when both bias and the initial proportion of agents with the superior opinion are low, the average consensus time is still polynomial in the network size.

Several directions remain open for theoretical investigation. One immediate problem is to theoretically establish the observed phase transitions for dense graphs. Here, it will be interesting to find out how `dense' a graph must be in order for it to exhibit the same phase transitions as in complete graphs. Similar questions remain open for sparse graphs such as $d$-regular expanders with constant $d$. Here, it will be of interest to find out how the threshold valued of the parameters depend on the average degree $d$ or the spectral properties of the expander. It is also worth analysing the exact behaviour below criticality. The numerical experiments conducted in this paper suggests that the consensus time is still polynomial below criticality (unlike in complete graphs where it is exponential), but it would be challenging to obtain exact bounds in this case.

\bibliographystyle{ieeetr}
\bibliography{opinion.bib}

\end{document}